\documentclass{amsart}
\usepackage{rlepsf,amssymb}

\oddsidemargin .25in \theoremstyle{plain}

\newtheorem{thm}{Theorem}

\newtheorem{cor}[thm]{Corollary}

\newtheorem{rem}[thm]{Remark}

\newtheorem{quest}[thm]{Question}

\newtheoremstyle{definition}{7pt plus6.3pt minus6.3pt}{7pt plus3pt minus3pt}%
{\rm}{}{\bf}{}{0.75em}{\thmname{#1}\thmnumber{ #2}\thmnote{\sl\stdspace#3}}
\theoremstyle{definition}\newtheorem{example}[thm]{Example}
\newtheorem{exercise}[thm]{\small Exercise}

\newcommand{\bbr}{\begin{rem}\em} 
\newcommand{\eer}{\end{rem}}

\newcommand{\bex}{\begin{example}} 
\newcommand{\eex}{\end{example}}

\newcommand{\bhw}{\begin{exercise}\small} 
\newcommand{\ehw}{\end{exercise}}

\newcommand{\be}{\begin{enumerate}}
\newcommand{\ee}{\end{enumerate}}

\def\tb{\operatorname{tb}}

\def\Z{\hbox{$\mathbb Z$} }

\def\R{\hbox{$\mathbb R$} }

\def\K{\mathcal{K}}
\def\L{\mathcal{L}}

\def\dfn#1{{\em #1}}

\begin{document}

\title{On contact surgery}

\author{John B. Etnyre}
\address{Georgia Institute of Technology, Atlanta, GA 30332-0160}
\email{etnyre@math.gatech.edu}
\urladdr{http://www.math.gatech.edu/\char126 etnyre}
\thanks{The author thanks Yasha Eliashberg for a helpful conversation during the preparation of this paper. 
Supported in part by NSF CAREER Grant (DMS--0239600) and FRG-0244663.}

\begin{abstract}
In this note we show that $+1$-contact surgery on distinct Legendrian knots frequently produces contactomorphic manifolds. We also give examples where this happens for $-1$-contact surgery. As an amusing corollary we find overtwisted contact structures  that contain a large number of distinct Legendrian knots with the same classical invariants and tight complements. 
\end{abstract}

\maketitle


\section{Introduction}

The purpose of this note is to study the behavior of contact structures under contact surgeries, with a particular emphasis on $+1$-contact surgeries. Legendrian surgeries
(that is $-1$-contact surgeries) have long been studied as a way to construct symplectically fillable, and hence tight, contact structures \cite{Eliashberg90a, Gompf98}. More recently, Ding and Geiges \cite{DingGeiges04} have  shown that any contact 3-manifold can be constructed from the standard tight contact structures on $S^3$ by a sequence of $\pm 1$-contact surgeries  ({\em cf} \cite{EtnyreHonda02a}). It was difficult to use this result at first because $+1$-contact surgeries frequently produce overtwisted contact structures. However, with the advent of the Heegaard-Floer invariants of contact structures \cite{OzsvathSzabo05a}, which behave very nicely under $+1$-contact surgeries, the Ding-Geiges result became an essential tool for constructing tight contact structures (especially tight but not symplectically fillable contact structures), \cite{LiscaStipsicz04}.
Despite these recent advances the nature of $\pm$-contact surgeries is still somewhat illusive. For example we have the following fundamental questions:\hfill\break
{\centerline {\em When does $\pm$-contact surgery preserve tightness?}}
and \hfill\break
\begin{center}{
\begin{minipage}[c]{3.2 in}
\centering{\em { When does $\pm$-contact surgery on distinct Legendrian knots give distinct contact structures?}}
 \end{minipage} \par}
\end{center}
Our main theorem begins to answer this question for $+1$-contact surgeries.
\begin{thm}\label{main}
Suppose $L_i, i=0,1,$ are Legendrian knots in a contact manifold $(M,\xi).$ If there exist a non-negative integer $n$ such that $S^n_\pm(L_0)$ is Legendrian isotopic to
$S^n_{\pm}(L_1)$ then $M_{+}(L_0)$ is contactomorphic to $M_{+}(L_1).$
\end{thm}
In this theorem, if $L$ is a Legendrian knot in a contact manifold $(M,\xi)$, then $S_\pm(L)$ denotes the $\pm$-stabilization of $L.$ In addition, $M_\pm(L)$ denotes the contact manifold that results from a $\pm 1$-contact surgery on $(M,\xi).$ 
Recall that given two Legendrian knots in the same topological knot type with the same Thurston-Bennequin invariants and rotation numbers then they become Legendrian isotopic after sufficiently many postive and negative stabilizations; however, there is no guarentee that the Legendrian knots will ever become the same after many postive stabilizations (or many negative stabilizations). In fact, most knot types in $S^3,$ for which we understand which stabilizations are necessary, only need stabilizations of one sign.  None the less, in Section~\ref{sec:ex} we discuss examples where Theorem~\ref{main} applies and examples where it does not.  
Moreover, in Corollary~\ref{trans}, we show that the hypothesis of Theorem~\ref{main} are met if the Legendrian knots are in a transversally simple knot type. 
In Theorem~\ref{gen} we generalize Thoerem~\ref{main}. With this generalization it is reasonable to ask if $+1$-contact surgery on Legendrian knots with the same classical invariants (knot type, Thurston-Bennequin invariant and rotation number) always produce contactomorphic manifolds.
We also see that this theorem provides the first examples of distinct Legendrian knots with the same classical invariants  on which $-1$-contact surgery yields contactomorphic manifolds. In particular, we prove the following result. 
\begin{cor}\label{firstcor}
Given any integer $n$ there is a tight contact manifold $(M,\xi)$ and $n$ distinct Legendrian knots $L_1,\ldots L_n$ in the same topological knot type with the same Thurston-Bennequin invariants and rotation numbers  such that all the contact manifolds obtained from $-1$-contact surgery on the $L_i$ are contactomorphic. 
\end{cor}
In Section~\ref{sec:ex}, we also ask several specific questions whose answers would significantly illuminate the nature of $\pm 1$-contact surgery. Finally our analysis of $+1$-contact surgery provides many examples of non-loose knots in overtwisted contact structures. Recall, a Legendrian knot in an overtwisted contact structure is called \dfn{loose} if its complement is also overtwisted. So far the only known example of non-loose knots are a family of unknots in one particular overtwisted contact structures on $S^3,$ \cite{Dymara??}. When the Thurston-Bennequin invariant and rotation number are fixed for these examples there is never more then two such examples (and maybe never more than one). However we can now show the following.
\begin{thm}\label{nonloose}
Given any integer $n$ there is an overtwisted contact manifold $(M,\xi)$ and distinct non-loose Legendrian knots $L_1,\ldots, L_n$ in the same knot type and having the same Thurston-Bennequin invariant and rotation number.
\end{thm}

\section{Contact surgery and notation}\label{background}

Throughout this note we use the standard results and terminology from convex surface theory. We refer the reader to \cite{EtnyreHonda01b, Honda00a} for all the necessary background.

Let $L$ be a Legendrian knot in a contact manifold $(M,\xi).$ Then $L$ has a \dfn{standard (contact) neighborhood} $N(L).$ This is a solid torus neighborhood of $L$ such 
that $\xi|_{N(L)}$ is tight and $\partial N(L)$ is a convex surface with two dividing curves each representing a longitude determined by the twisting of the contact planes along $L$ (that is the framing of $L$ coming form $\xi$ and from a dividing curve agree). This neighborhood is uniquely determined by $L$ up to isotopy through such neighborhoods. Identify $\partial N(L)$ with $\R^2/\Z^2$ so that $\pm(1,0)^T$ corresponds to the boundary of the meridional disk and $\pm (0,1)^T$ corresponds to the curve determined by the dividing curves. We call these standard coordinates on the boundary of a standard neighborhood. Identifying $\partial {M\setminus N(L)}$ with
$-\partial N(L)$ we can define the map $\phi_\pm: \partial (D^2\times S^1)\to \partial ({M\setminus N(L)})$ by 
\[\phi(x,y)= \begin{pmatrix}1& 0\\ \pm 1 & 1\end{pmatrix}\begin{pmatrix}x\\ y\end{pmatrix}.\] 
Let $M_\pm(L)$
be the manifold obtained by gluing $D^2\times S^1$ to $M\setminus N(L)$ using this map. The contact structure $\xi$ induces a contact structure on $M\setminus N(L)$ and 
the two dividing curves on $\partial (M\setminus N(L)),$ as seen on $\partial (D^2\times S^1),$ represent $(\mp 1, 1)$ curves. Thus, according to \cite{Honda00a}, there is a unique tight contact structure on $D^2\times S^1$ having convex boundary with these dividing curves. So we may extend $\xi|_{M\setminus N(L)}$ to a contact structure 
$\xi_\pm$ on $M_\pm.$ The contact manifold $(M_\pm,\xi_\pm)$ is said to be obtained from $(M,\xi)$ by \dfn{$\pm 1$-contact surgery on $L$}. The term \dfn{Legendrian
surgery} typically refers to $-1$-contact surgery. 

\section{Proof of Theorem~{\protect \ref{main}}}

We are now ready to prove our main theorem.
\begin{proof}[Proof of Theorem~\protect \ref{main}]
Assume the hypothesis of the theorem hold with positive stabilizations (the proof for negative stabilziations is analogous). 
Let $M_i=M\setminus N(L_i).$ (Recall, when $L$ is a Legendrian knot, $N(L)$ will refer to a standard neighborhood of $L.$) Moreover, let $M'_i=M\setminus N(S^n_+(L_i)).$ Recall that we may assume that $N(S^n_+(L_i)) \subset N(L_i),$ see \cite{EtnyreHonda01b}.
Since $S_+^n(L_0)$ is Legendrian isotopic to $S_+^n(L_1)$ we know $M_0'$ is contactomorphic to $M_1'.$ Moreover, $A_i=M_i'\setminus M_i$ is diffeomorphic to $T^2\times [0,1]$ and the contact structures on $A_0$ and $A_1$ are contactomorphic, universally tight, positive continued fractions blocks (see \cite{Honda00a} for terminology). 
More specifically if $A_i\cap M_i=T^2\times\{1\}$ then the dividing curves on $T^2\times\{1\}$ are in the class $\pm(0,1)^T,$ the dividing curves on $T^2\times \{0\}$
are in the class $\pm(m,1)^T$ (here we are using the standard coordinates on $\partial N(L)$ as described in Section~\ref{background} to identify $T^2$ with $\R^2/\Z^2$) and all bypasses on a vertical annulus in $A_i$ are positive. 

Let $(W_i,\xi_i)$ be the result of $+1$-contact surgery on $(M,\xi)$ along $L_i$ and let $L_i'$ be the Legendrian knot in the core of the surgery torus in $(W_i,\xi_i).$
So $M_i=W_i\setminus N(L_i').$ Note that in $N(L'_i)$ there is a torus $T_i$ that splits $N(L_i')$ into $A'_i$ and $S_i$ where $A_i'=T^2\times [0,1]$ with $T^2\times{1}=N(L_i')\cap M_i$ and $T^2\times\{0\}=T_i,$ $S_i$ is a solid torus, the dividing curves on $T_i$ are in the class $\pm(m, m+1)$ and $A_i'$ is a universally tight positive continued fractions block. (We are using standard coordinates on $\partial N(L')$ as described above to identify $T^2$ with $\R^2/\Z^2.$) Thus, as contact manifolds,
we have
\[W_i=M_i\cup N(L_i')=M_i\cup A_i'\cup S_i,\]
but  notice that $M_i\cup A_i'$ is contactomorphic to $M_i\cup A_i$ (since $A_i$ and $A_i'$ are the same continued fractions blocks expressed in different coordinates)
so
\[W_i=M_i\cup A_i\cup S_i=M_i'\cup S_i. \]
Now we know $M_0'$ is contactomorphic to $M_1'$ and so we are left to see this contactomorphism can be extended over the $S_i.$ To this end notice that there are exactly two tight contact structures on a solid torus with convex boundary and dividing curves as given by the $S_i$'s. These two contact structures are distinguished by the sign of a bypass on the meridional disk, see \cite{Honda00a}. Since $N(L_i')$ is tight the sign of the bypass in $S_i$ must agree with the sign of the bypasses (along a horizontal annulus) in $A_i'.$ Thus the contact structures on $S_0$ and $S_1$ are contactomorphic and any diffeomorphism of the boundary preserving the characteristic 
foliation can be extended to a contactomorphism. Thus we may extend our contactomorphism from $M'_0$ to $M'_1,$ to a contactomorphism from $W_0$ to $W_1.$
\end{proof}
\bbr
Note that a continued fraction block containing basic slices of the same sign can be embedded into solid tori with various slopes. It turns out that if the signs in the continued fractions block are different then this seriously limits the way in which the continued fractions block can be embedded, \cite{Honda00a}. This is the reason why we cannot generalize Theorem~\ref{main} to allow for more general stabilizations of the Legendrian knots. 
\eer

\section{Examples, corollaries and further discussion}\label{sec:ex}

\begin{figure}[ht]
  \relabelbox 
  \small{\epsfxsize=2in\centerline{\epsfbox{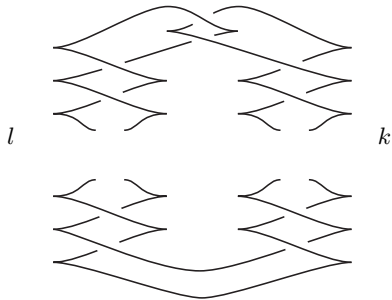}}}
  \relabel {l}{$l$}
  \relabel {k}{$k$}
  \endrelabelbox
  \caption{The Legendrian knot $E_n(k,l)$. There are $l$ crossings on the right
    hand side and $k$ on the left.}
  \label{fig:gench}
\end{figure}

\begin{figure}[ht]
  \relabelbox 
  \small{\epsfxsize=3in \centerline{\epsfbox{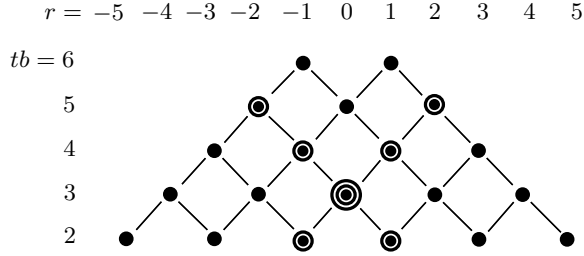}}}
  \relabel {- 5}{$-5$}
  \relabel {- 4}{$-4$}
  \relabel {- 3}{$-3$}
  \relabel {- 2}{$-2$}
  \relabel {- 1}{$-1$}
  \relabel {0}{$0$}
  \relabel {1}{$1$}
  \relabel {2}{$2$}
  \relabel {3}{$3$}
  \relabel {4}{$4$}
  \relabel {5}{$5$}
  \relabel {tb = 6}{$tb=6$}
  \relabel {a}{$5$}
  \relabel {b}{$4$}
  \relabel {c}{$3$}
  \relabel {d}{$2$}
  \relabel {r}{$r=$}
  \endrelabelbox
 \protect{ \caption{Classification of Legendrian $(2,3)$-cables of $(2,3)$-torus 
    knots.  Concentric circles indicate multiplicities, {\em ie}, the number of 
    distinct isotopy classes with a given $r$ and $\tb$.\hfill}   }
  \label{iterate}
\end{figure}

In this section we explore applications of Theorem~\ref{main} to $\pm 1$-contact surgeries. We begin with a much studied set of examples.
\bex\label{firstex}
Given an integer $n$ larger than three and two positive integers $l$ and $k$ such that $n=l+k$ we have the Legendrian knot $K_n(k,l)$ shown in Figure~\ref{fig:gench}.
For a fixed $n$ we get a single topological knot type, but using Legendrian contact homology \cite{EliashbergGiventalHofer00} one can show \cite{Chekanov02, EpsteinFuchsMeyer01} 
that $K_n(k,l)$ is Legendrian isotopic to $K_n(k',l')$ if and only if the unordered pair $\{k,l\}$ is the same as the unordered pair $\{k',l'\}.$ 
It is natural to ask if the result of $\pm 1$-contact surgery on the $K_n(l,k)$'s yield distinct contact structures. In particular, one could then use invariants of the contact
structures, like the Heegaard-Floer contact invariants \cite{OzsvathSzabo05a}, on the surgered manifolds to try to distinguish the knots. 
In \cite{LiscaStipsicz06?} it was shown that the Heegaard-Floer invariants $\hat{c}(S^3_+(K_n(k,l))$ were the same for all $k$ and $l$ (with fixed $n$ of course). Theorem~\ref{main} shows
why this must have been the case. Specifically, it is easy to verify (or see \cite{EpsteinFuchsMeyer01}) that $K_n(k,l)$ and $K_n(k-1, l+1)$ are isotopic after a single positive stabilization if $l$ is even and after a single negative stabilization if $l$ is odd. Thus Theorem~\ref{main} implies the following result. 
\begin{cor}
The contact manifold $S^3_+(K_n(k,l))$ is independent of $k$ and $l$.
\end{cor}
We are unable to say anything about the contactomorphism type of $S^3_-(K_n(k,l)),$ but in \cite{LiscaStipsicz06?} it is shown that the Heegaard-Floer invariants of the $S^3_-(K_n(k,l))$ are the same when $k$ is odd. The following question is still open.
\begin{quest}
Are the contact manifold $S^3_-(K_n(k,l)),$ for $k=1,\ldots, \lfloor\frac n2\rfloor$, distinct?
\end{quest}
\eex
It is clear from the last example that the contact manifold obtained by $+1$-contact surgery on a Legendrian knot
does not provide a very strong invariant of the Legendrian knot. We can illustrate this further with the following result. 
\begin{cor}\label{trans}
If $\mathcal{K}$ is a topological knot type that is transversely simple (recall this means transverse knots in this knot type are determined by their self-linking number) and 
$L_0$ and $L_1$ are Legendrian knots in this knot type with the same Thurston-Bennequin invariant and rotation number then $S^3_+(L_0)$ is contactomorphic to 
$S^3_+(L_1).$
\end{cor}
\begin{proof}
Since $\mathcal{K}$ is transversally simple the transverse push off of $L_0$ and $L_1$ are transversely isotopic and hence the Legendrian knots themselves are 
Legendrian isotopic after a suitable number of negative stabilizations \cite{EtnyreHonda01b}. Thus the result follows form Theorem~\ref{main}.
\end{proof}
Most knot types for which we have a classification are transversely simple, currently the only two non-transversely simple knot types come from 3-braids \cite{BirmanMenasco06II} and a cable 
of a torus knot \cite{EtnyreHonda05}. As there is no classification of Legendrian knots in the former knot type we restrict attention to the later knot type
\bex
Let $\mathcal{K}$ be the knot type of the $(2,3)$-cable of the $(2,3)$-torus knots. The classification of knots in $\mathcal{K}$ form \cite{EtnyreHonda05} is indicated in Figure~\ref{iterate}. 
Here the dots represent Legendrian knots with a given Thurston-Bennequin invariant and rotation number, a dot with a circle around it represents two distinct Legendrian knots with the given invariants and the edges represent positive and negative stabilizations. 
Note that given any two Legendrian knots with the same invariants then after either positive or negative stabilization they will be isotopic, thus the result of $+1$-contact
surgery on these knots will always produce contactomorphic manifolds. In particular, consider the two Legendrian knots $L$ and $L'$ with $tb= 5$ and $r=2.$ In \cite{LiscaStipsicz06?}, Lisca and Stipsicz showed that the Heegaard-Floer invariants of $S^3_+(L)$ and $S^3_+(L')$ both vanish. Since one of the Legendrian knots, say $L,$ is a stabilization it is easy to see that $S^3_+(L)$ is overtwisted, which in turn implies the vanishing of the Heegaard-Floer invariants. It is not, {\em a priori}, clear that $S^3_+(L')$ is overtwisted, but according to Theorem~\ref{main}, it is.  
\eex

It is well known that any two Legendrian knots in the same knot type with the same Thurston-Bennequin invariant and rotation number are Legendrian isotopic after they have been positively and negatively stabilized sufficiently many times. Considering the examples above and Theorem~\ref{main} one might now ask if any two knots with the same Thurston-Bennequin invariant and rotation number become isotopic after positive or negative (but not both) stabilizations? This turns out not to be the case.
\bex
Again consider the $(2,3)$-cable of the $(2,3)$-torus knot. Let $L$ and $L'$ be the two 
Legendrian knots in this not type with $tb=5$ and $r=2$ and let $K$ and $K'$ be the two Legendrian knots with $tb=5$ and $r=-2.$ Moreover, assume that $L'$ and $K'$ are the Legendrian knots that are not stabilizations of other Legendrian knots. Finally set $E=L\#K$ and $E'=L'\#K'.$ The knots $E$ and $E'$ have $tb=11$ and $r=0.$ 
In \cite{EtnyreHonda03} connect sums of Legendrian knots were studied and it was shown that if $\mathcal{L}(\mathcal{K})$ represent Legendrian knots in the knot type $\mathcal{K}$ and $\K=\K_1\#\ldots \#\K_n$ is a topological connected sum knot type in $S^3,$
with $\K_i$ prime, then the map
\[\frac{\L(\K_1)\times\ldots \times\L(\K_n)}{\sim} \to \L(\K_1\# \ldots\#\K_n)\]
given by connect sum is a one to one correspondence where $\sim$ is generated by
\begin{enumerate}
\item $(\ldots, S_\pm(L_i),\ldots, L_j, \ldots)\sim(\ldots, L_i,\ldots, S_\pm(L_j), \ldots)$ and
\item $(L_1,\ldots, L_n)\sim(L_{\sigma(1)},\ldots, L_{\sigma(n)})$ where $\sigma$ is permutation
        of $1, \ldots, n$ such that $\K_i=\K_{\sigma(i)}.$ 
\end{enumerate}
From this it is clear that $E$ and $E'$ remain distinct Legendrian knots under positive or negative stabilization (thought of course they become the same after a positive and a negative stabilization).  Thus from our theorem it is not immediately clear if $S^3_+(E)$ is contactomorphic to $S^3_+(E')$. But consider $E''=L\# K'$ and notice that
\[S_+(E'')=(S_+(L))\# K'=(S_+(L'))\#K'=S_+(E').\]
Similarly $S_-(E'')=S_-(E).$ Thus $S^3_+(E)=S^3_+(E'')=S^3_+(E').$
\eex
This last example shows how to generalize Thoerem~\ref{main} a little. 
\begin{thm}\label{gen}
Suppose $L$ and $L'$  are Legendrian knots in a contact manifold $(M,\xi).$ If there exist Legendrian knots $L_1,\ldots , L_k,$ non-negative integers $n_1,\ldots, n_{k-1}$ and signs $\epsilon_1,\ldots, \epsilon_{k-1}$ such that $L_1=L, L_k=L'$ and $S_{\epsilon_i}^{n_i}(L_i)$ is Legendrian isotopic to $S_{\epsilon_i}^{n_i}(L_{i+1}),$ for $i=1,\ldots, k-1,$ then $M_+(L)$ is contactomorphic to $M_+(L').$
\end{thm}
This prompts the following question.
\begin{quest}
Are two Legendrian knots with the same classical invariants (knot type, Thurston-Bennequin invariant and rotation number) always related by a set of knots as in the hypothesis of Theorem~\ref{gen}?
\end{quest}
Of course a positive answer to this question would imply a positive answer to the following fundamental question.
\begin{quest}
Does $+1$-contact surgery on two Legendrian knots with the same classical invariants always produce contactomorphic contact manifold?
\end{quest}

Our discussion so far has indicated that $+1$-contact surgery is quite a bad invariant of Legendrian knots, but we have gained no insight into $-1$-contact surgeries. There is very little known here, but it is much more likely that $-1$-contact surgeries on distinct Legendrian knots yield distinct manifolds. That said it is not always the case that the manifolds are distinct. Recall, Corollary~\ref{firstcor} says: 
Given any integer $n$ there is a tight contact manifold $(M,\xi)$ and $n$ distinct Legendrian knots $L_1,\ldots L_n$ in the same topological knot type with the same invariants
such that all the contact manifolds obtained from $-1$-contact surgery on the $L_i$ are contactomorphic. 
\begin{proof}[Proof of Corollary~\ref{firstcor}]
Consider again the twist knots $K_n(k,l)$ from Example~\ref{firstex}. When $n$ is even $K_n(k,l)$ has Thurston-Bennequin invariant 1 and rotation number 0. It is also known that $K_n(k,l)$ has slice genus $g_s(K_n(k,l))=1.$ Thus $tb(K_n(k,l))=2g_s(K_n(k,l))-1$ and a theorem of Lisca and Stipsicz \cite{LiscaStipsicz04} implies that $S^3_+(K_n(k,l))$ is tight for all even $n.$ For a fixed $n$ the above mentioned result of \cite{EpsteinFuchsMeyer01} implies there are $\frac n2$ distinct Legendrian knots $K_n(k,l).$ In fact they are $K_n(k,l)$ for $k=1,\ldots, \frac n2.$ Now let $M_k$ be the result of $+1$-contact surgery on $K_n(k,l)$ for $k=1,\ldots, \frac n2.$ From our 
discussion above we see that all the $M_k$'s are contactomorphic, so we denote their common contactomorphism class by $M.$ However, the core of each surgery torus 
gives a Legendrian knot $K_k$ in $M.$ Of course, $M\setminus N(K_k)$ is contactomorphic to $S^3\setminus K_n(k,l)$ and hence for distinct $k$ between 1 and $\frac n2$
the contact manifolds $M\setminus N(K_k)$ are not contactomorphic. Thus the Legendrian knots $K_k$ in $M$ are not Legendrian isotopic and we have $\frac n2$ distinct Legendrian knots in the same topological knot type with the same Thurston-Bennequin invariant and rotation number. Note that a $-1$-contact surgery on the $K_k$ will undo the the original surgery, thus returning us to $S^3$ with the standard contact structure. Thus we have found $\frac n2$ distinct Legendrian knots on which $-1$-contact surgery yields the same contact manifold. 
\end{proof}
It is much more difficult to find such examples in $S^3$. So we merely ask the following open question.
\begin{quest}
Does $-1$-contact surgery on distinct Legendrian knots in the standard tight contact structure on $S^3$ always produce distinct contact manifolds?
\end{quest}

\section{Non-loose knots}
We are now ready to produce families of non-loose Legendrian knots with the same classical invariants in overtwisted contact structures.
\begin{proof}[Proof of Theorem~\ref{nonloose}]
In \cite{EtnyreHonda03} the following theorem was proven.
\begin{thm}\label{thm:examples} 	
Given two positive integers $m$ and $n$, there are
distinct Legendrian knots $L_1,\dots, L_n$ in the same topological knot type which 
have the same Thurston-Bennequin invariant and rotation number, and remain 
distinct even after $m$ stabilizations (of any type). 
\end{thm}
The examples in this theorem were constructed by taking the connected sum of various negative torus knots in $S^3.$ Theorem~1.1 in \cite{LiscaStipsicz06?} implies that $+1$-contact surgery on all these examples produce overtwisted contact structures. Since all the $L_i$ have the same invariants, it is easy to see the resulting contact structures are all homotopic. Thus $S^3_+(L_i)$ is independent of $i.$ Denote the resulting overtwisted contact manifold by $M.$ The core of the surgery torus associated to surgery on $L_i$ gives a Legendrian knot $L_i'$ in $M.$ All the $L_i'$'s are topologically isotopic in $M$ and have the same classical invariants. Moreover, the complement of $L_i'$ in $M$ is contactomorphic to the complement of the $L_i$ in $S^3$ and thus the $L_i'$ are non-loose knots. Moreover the contact structure on $M\setminus L_i$ is distinct for each $i,$ since the $L_i$ in $S^3$ are distinct knots and two Legendrian knots in $S^3$ are Legendrian isotopic if and only if their complements are contactomorphic, \cite{Etnyre05}. Thus the $L_i'$ are the desired examples. 
\end{proof}

\def\cprime{$'$} \def\cprime{$'$}

\end{document}